\documentclass{elsart}
\usepackage{natbib}
\usepackage{amsfonts}

\begin{document}

\begin{frontmatter}



\title{Graphs with restricted valency and matching number}


\author{Niranjan Balachandran},
\address{The Ohio State University, Columbus, Ohio, USA }
\ead{niranj@math.ohio-state.edu}

\author{Niraj Khare\corauthref{cor}}
\corauth[cor]{Corresponding author}
\ead{nirajkhare@math.ohio-state.edu}

\address{The Ohio State University, Columbus, Ohio, USA }

\begin{abstract}
% Text of abstract
Consider the family of all finite graphs with maximum degree $\Delta(G)<d$ and matching number $\nu(G)<m$. In this paper we give a new proof to obtain the exact upper bound for the number of edges in such graphs and also characterize all the cases when the maximal graph is unique. We also provide a new proof of Gallai's lemma concerning factor critical graphs.
\end{abstract}


\end{frontmatter}

% main text
\section{Introduction}
Let $\mathcal{F}(d,m)$ denote the set of all finite maximal simple graphs that satisfy $\Delta(G)<d,\nu(G)<m$, where $\Delta(G)$ denotes the maximum degree among all the vertices of $G$ and $\nu(G)$ denotes the maximum matching size in $G$. Here, maximality is with respect to the edge set, i.e., if a graph $G$, a member of $\mathcal{F}(d,m)$, is a subgraph of $G'$ with $|E(G')|>|E(G)|$ then either $\Delta(G')\geq d$ or $\nu(G')\geq m$. In particular, when $d=m=s$, this set consists of all those finite maximal graphs with both degree and matching size less than $s$. In other words, these are set systems with uniform block size $2$ containing no sunflower with $s$ or more petals. A \textit{sunflower} with $s$ petals is a collection of sets $A_1,A_2\ldots, A_s$ and a set $X\textrm{(possibly empty)}$ such that $A_i\cap A_j=X$ whenever $i\neq j$. The set $X$ is called the core of the sunflower.\\ \\
 It is a well known result (due to Erd\H{o}s-Rado\cite{ER}) that a uniform set system with block size $k$ of size greater than $k!(s-1)^k$ admits a sunflower with $s$ petals(for a proof see \cite{BF}, for instance). Other bounds that ensure the existence of a sunflower with $s$ petals are known in the case of $s=3$ with block size $k$( Kostochka\cite{K}, for instance), but the general case seems quite far away.\\
 Our instance of this problem deals with $k=2$ (the graph case) and with $s$, arbitrary. The Erd\H {o}s - Rado bound in this case ($2(s-1)^2$) is trivially obtained and our result yields the exact bound on the number of edges. This was first achieved in Sauer \emph{et al}\cite{AH}.\\ 
 \hspace{2cm} The exact bound for the general case(arbitrary $d,m$) was first obtained in \cite{CH}. Our proof is simpler, self-contained and  `structural' as opposed to Chv\'atal \emph{et al}\cite{CH}. Moreover, our method enables us to give a simple characterization of all the cases where the subfamily of $\mathcal{F}(d,m)$ $-$graphs having no isolated vertices$-$ attaining the bound  is a singleton, i.e., a unique graph, up to isomorphism, attains the bound. We also present a proof of Gallai's lemma enroute as an application of the main theorems proved in section 3. We use the method of augmenting paths in graphs to study maximal matchings. 
\section{The Problem}
\noindent Throughout this paper, we shall denote by $e_{xy}$, the edge $\{x,y\}$ in a graph $G$ with $x,y$ being vertices of $G$. By a \textit{walk} in a graph $G$ with vertex set $V$ and edge set $E$, we shall mean an ordered sequence $(x_1,x_2,\ldots,x_n)$ where $x_i$ are vertices of $G$, and $e_{x_ix_{i+1}} \in E$. If in addition we also have $x_i\neq x_j$ for $i\neq j$, the walk is a \textit{path} of $G$.
We call a vertex $v$ \textit{unsaturated} relative to a matching $\mathcal{M}$ if $v$ is not covered by $\mathcal{M}$.\\
\begin{defn}: Let $d\geq 2, m\geq 2$ be integers. Let  $\mathcal {F} (d,m) $ be the family of graphs such that $G \in \mathcal{F} (d,m)\Leftrightarrow $\\
(i) $\Delta (G)< d $ and $\nu (G) < m $.\\
(ii) $E(G)$ is maximal with property (i).
(i.e., if $G'$ is a graph s.t. $G$ is a subgraph of $G'$ and $|E(G')| > |E(G)|$ then either $\Delta( G')\ge d$ or $\nu(G')\geq m$).\end{defn}
Note that $G \in\mathcal{F} (d,m) $\ $\Rightarrow$ \ $ |E(G)|\le  2(m-1)(d-1) - (m-1)$. This is immediate since the unsaturated vertices of any such graph form an independent set and further, the set of vertices of a maximal matching $\mathcal{M}$ gives us  a vertex cover for the edges of $G$. \\ \\
We define for any graph $G$, \\
$V_{\ge  1} (G) : = \{v \in V(G) |  d_{G} (v) \ge 1\}, V_{0} (G)      : =  \{v \in V(G) |  d_{G} (v) = 0\}$
and let $e(d,m)\ :=\displaystyle\textrm{Sup}_{G\in\mathcal{F} (d,m)} |E(G)|$. We aim at obtaining a precise upper bound for $e(d,m)$. In particular, for the case where $d=m=s$, for $s \ge  4$, we are dealing with the Erd\"os-Rado problem of sunflowers with $s$ petals for the case $k=2$.\\
It is easy to see that $e(2,2) = 1$ and ${\mathcal {F}}  (2,2)$ has a unique graph, with no isolated vertices, $K_2$ as a member. Likewise, $e(3,3) = 6$ and ${\mathcal {F}}  (3,3)$ has a unique member, with no isolated vertices,  consisting of two disjoint $K_3's$.\\
We make the following very simple observation before moving on.
\begin{prop}: If $G \in\mathcal{F} (d,m)$ then $\nu (G) = m-1$.\end{prop}
\textbf{Proof}:  If not then $\nu (G) < m-1$ since $\nu (G) \le m-1$.\\
Define a new graph $G'$:\\
$V(G') = V(G) \cup  \{1, 2 \}$ where $1,2 \notin  V(G), E(G') = E(G) \cup  e_{12}$.\\
 Note that $\Delta (G') = \max \{ \Delta (G), 1\}  < d$ and $ \nu (G')  = \nu (G) + 1 < m$. This contradicts the assumption that $G \in \mathcal{F} (d,m)$ since $|E(G')| > |E(G)|$ . \hspace{\stretch{1}} $\square$  
\section{Preliminaries}
Let $G$ be a simple graph and $\mathcal {M} $ be a matching of $G$. Let $|\mathcal {P} |$ denote the number of edges in the path $\mathcal{P}$ of $G$.
\begin{defn}: A path $\mathcal {P}$ of G is called a star path of $G$ relative to $\mathcal {M}$ if and only if $\mathcal {P}$ satisfies the following:\\
(i) $|\mathcal {P}|\ \equiv \ 0\pmod 2$,\\
(ii) $\mathcal {P}$ has alternating non-matching and matching edges $w.r.t\  \mathcal {M}$,\\
(iii) $\mathcal{P}$ starts with an unsaturated (relative to $\mathcal {M}$) vertex.\end{defn}

Star paths of $G$ relative to $\mathcal{M}$ of length zero are simply the unsaturated vertices of $G$.

\begin{defn}: A $v \in V(G)$ is called a \textit{star vertex}, relative to $\mathcal {M}$ if a star path of $G$ terminates at $v$.\end{defn}
\begin{defn}: $Star(G,\mathcal {M}) := \{v \in \ V(G)| v \textrm{\ is a star vertex of\ } G\}$.\end{defn}
 If $v \in  V(G)$ and $v$ is unsaturated relative to $\mathcal{M}$ then $v \in  Star(G,\mathcal{M})$. This is immediate since there exists a star path of length zero terminating at $v$.

We shall now set up some more notation. For any path $\mathcal {Q} =( x_1, x_2, x_3, ...,x_n) $ of G let $\overline {\mathcal {Q}}  :=( x_n,...,x_3,x_2,x_1)$ denote the reverse path. For paths $\mathcal{P}=(x_1,x_2\ldots,x_n),\mathcal{Q}=(y_1,y_2,\ldots,y_m)$ such that $x_n = y_1$, we denote by $\mathcal{P}\star \mathcal{Q}$, the walk, $(x_1,\ldots,x_n,y_2,\ldots,y_m)$  resulting as a concatenation of the paths $\mathcal{P}$ and $\mathcal{Q}$.
\begin{prop}: Let $G$ be a given graph with a maximal matching $\mathcal{M}$. Let $v \in  V(G)\setminus Star(G, \mathcal {M})$ and $u \in  V(G)$. Define $G'$ as follows: $V(G') := V(G)$ , $E(G') := E(G) \cup  \{e_{uv}\}$. Then $\nu(G')\ =\ \nu(G)$. Further, $\mathcal{M}$ is a maximal matching of $G'$ and $v\in  V(G')\setminus  Star(G', \mathcal {M})$.\end{prop}
\textbf{Proof}: If $e_{uv} \in E(G)$ then there is nothing to prove. So, we assume $e_{uv} \notin E(G)$.\\ 
We show that  $\mathcal{M}$ is a maximal matching of $G'$.\\
Suppose not. Then by the theory of augmenting paths(see \cite{PS} for instance), there is an $\mathcal{M}$ augmenting path, $\mathcal {P}$ which necessarily contains the edge $e_{uv}$ else $\mathcal{P}$ would be a path in $G$, which is not possible as $\mathcal{M}$ is maximal in $G$. Let  $\mathcal{P}=( x_1, x_2, ...,x_n,u,v,w,y_m,...,y_2,y_1)$ in $G'$ be the augmenting path. Note that since unsaturated ($w.r.t\ \mathcal{M}$) vertices of a graph are always elements of $Star(G,\mathcal{M})$, 
and $v \notin Star(G, \mathcal{M})$, it follows that $v$ is saturated by $\mathcal{M}$. Hence an augmenting path for $\mathcal{M}$ cannot terminate at $v$ in $G'$. This justifies the writing of $w,y_m\ldots,y_1$ in the augmenting path.
Also note that $e_{uv} \notin \mathcal{M}$ since $e_{uv}\notin E(G)$ by assumption. Consequently $e_{wv}\in\mathcal{M}$.\\
Consider $\mathcal{Q}\ :=\ (y_1,y_2,...,y_m,w,v)$. This walk is in fact a star path in $G$ terminating at $v$. But that is a contradiction since $v\ \in \ V(G)\setminus Star(G, \mathcal {M})$. \\ 
We next show that $v\in V(G')\setminus Star(G',\mathcal {M})$. Suppose not, then there exists a star path $\mathcal {Q}$ in $G' $ terminating at $v$ containing the $\mathcal {M}$-edge incident at $v$ as the terminal edge of the path. 
Hence the non-matching edge, $e_{uv}$ does not belong to $\mathcal{Q}$. But then $\mathcal {Q}$ lies in $G$ and hence $v$ is a star vertex of $G$ as well, which by the assumption is a contradiction. 
\hspace{\stretch{1}} $\square$\\  \\
Before we prove the next lemma, we make a remark: If $G \in  \mathcal {F}(d,m)$, $\mathcal{M}$ be a matching of $G$ and $v \in  V(G)\setminus Star(G, \mathcal {M})$ then $deg(v)= d-1$. Indeed, suppose there exists $ v \in  V(G)\setminus Star(G, \mathcal {M})$ such that $deg(v) < d-2$. Assume, without loss of generality, that there exists $u \in V_{0}(G)$. Define a new graph $G'$, $V(G') = V(G)$ and $E(G') = E(G) \cup  e_{uv}$.
Therefore $\nu (G') = \nu (G)$ by Proposition 6 as $v \in  V(G)\setminus Star(G, \mathcal {M})$.
But since $|E(G')| > |E(G)|$, we have a contradiction. \\
Next, we prove a lemma which will lead us to our main result-- a connected graph, all of whose vertices are star vertices relative to a maximum matching, has a unique vertex which is not covered by the matching-- the theorem 8.
\section{The Main Result}
\begin{lem}: Let $G$ be a simple graph and $\mathcal{M}$, a matching of $G$.
Let \begin{displaymath}\mathcal{P}_1:=(x_1,x_2,\ldots,x_n),\  \mathcal{P}_2:=(y_1,y_2,\ldots,y_m)\end{displaymath} be star paths in $G$. If $V(\mathcal{P})_1 \cap V(\mathcal{P}_2) \neq \emptyset$ then either there exists a star path from $x_1$ to $y_m$ or there exists an $\mathcal{M}$-augmenting path from $x_1$ to $y_1$.\end{lem}
\textbf{Proof}: Let
 \begin{displaymath}  i=\min\{k | x_k\in\mathcal{P}_1\cap\mathcal{P}_2\}  \end{displaymath}
 Without loss of generality $i>1$. If not, then $x_1\in\mathcal{P}_1\cap\mathcal{P}_2$. As $x_1$ is not covered by $\mathcal {M}$ it follows that $x_1=y_1$. This is so since $y_1$ is the only vertex in $\mathcal{P}_2$ which is not covered by $\mathcal {M}$. So there is a star path(in this case $\mathcal{P}_2$) from $x_1$ to $y_m$.\\
 So $x_i=y_j$ for some fixed $j \in \{ 2,3, \ldots,m\}$. Note that $y_1\notin\mathcal{P}_1\cap\mathcal{P}_2$  as $i>1$ and $y_1$ is not covered by $\mathcal {M}$. Furthermore if $j=m$, then by the definition of $i$ and the fact that $\mathcal{P}_1$, $\mathcal{P}_2$ are star paths, we have an $\mathcal {M}$ augmenting path $(x_1, \ldots,x_i) \star \overline { \mathcal{P}_2}$ from $x_1$ to $y_1$. So let  $j \in \{ 2,3, \ldots,m-1\}$\\
 Now we observe that $e_{x_{i-1}x_i}\notin\mathcal{M}$. Suppose not. Since $x_{i-1}\notin\mathcal{P}_1\cap\mathcal{P}_2$, we have $x_{i-1}\neq y_{j-1}$ and also $x_{i-1}\neq y_{j+1}$. Thus, $e_{y_{j-1}y_j}\notin\mathcal{M}$ and $e_{y_jy_{j+1}}\notin\mathcal{M}$ since there is atmost one matching edge incident at $x_i=y_j$. But that is a contradiction to the fact that $\mathcal{P}_2$ is also alternating path relative to $\mathcal {M}$.\\  
 As $e_{x_{i-1}x_i}\notin\mathcal{M}$, we have $e_{y_{j-1}y_j}\in\mathcal{M}$ or $e_{y_jy_{j+1}}\in\mathcal{M}$.\\ \\
 Case 1: $e_{y_{j-1}y_j}\in\mathcal{M}$.\\
 In this case, consider the sequence of vertices $\mathcal{P}:=(x_1,x_2\ldots,x_i,y_{j-1},\ldots, y_1)$. By the definition of $i$, $y_k\notin\{x_1,x_2,\ldots,x_{i-1}\}$ for all $1\leq k\leq m$. Using the fact that $\mathcal{P}_1,\mathcal{P}_2$ are star paths, it now easily follows that $\mathcal{P}$ is an augmenting path for $\mathcal{M}$.\\ \\
 Case 2: $e_{y_jy_{j+1}}\in \mathcal{M}$.\\
  In this case, let $\mathcal{P}:=(x_1,x_2,\ldots,x_i,y_{j+1},\ldots,y_m)$. As before, by the definition of $i$, $y_k \notin \{x_1,x_2\ldots,x_{i-1}\}$ for $1\leq k\leq m$ and the fact that $\mathcal{P}_1,\mathcal{P}_2$ are star paths implies that $\mathcal{P}$ is also a path and in fact is a star path from $x_1$ to $y_m$ as desired. \hspace{\stretch{1}} $\square$
\begin{thm}: Let $G$ be a simple graph and $\mathcal{M}$ a maximum matching of $G$. If there is a connected component $\mathcal{C}$ of $G$ such that $V(\mathcal{C})\subseteq Star(G,\mathcal{M})$ then there exists \underline {exactly one} \textit{unsaturated} ( relative to $\ \mathcal{M}$) vertex in $\mathcal{C}$.\end{thm}
\textbf{Proof}: Note that $\mathcal{C}$ has at least one unsaturated vertex in it. Indeed, since every vertex is a member of $Star(G,\mathcal{M})$, there is at least one unsaturated vertex in $\mathcal {C}$ for star paths to originate from. We now prove that $\mathcal{C}$ has at most one unsaturated vertex.\\ Suppose $v\in V(\mathcal{C})$ is an $\mathcal{M}$-unsaturated vertex.\\ Define $T_v:=\{y\in V(\mathcal{C})| \textrm{\ there is a star path from\ } v \textrm{\ to\ } y\}$. We claim that $T_v=V(\mathcal{C})$. Note that $v \in T_v$, so $T_v \neq \emptyset$. Let $N(u)$ denote the set of neighbors of a vertex $u$. \\ 
 Suppose $T_v \neq V(\mathcal {C})$. Then there exists a vertex $u\in V(\mathcal{C})\setminus T_v$ such that $N(u)\cap T_v\neq\emptyset$. Let $w\in N(u)\cap T_v$. By the definition of $T_v$, there exists a star path $\mathcal{P}_1:=(v,x_1,\ldots,,x_n,w)$.\\ 
 As $u\in V(\mathcal{C})\subseteq Star(G,\mathcal{M})$, there exists a star path $\mathcal{P}_2:=(y,y_1,\ldots,y_m,u)$. Note if $\mathcal{P}_1\cap \mathcal{P}_2=\emptyset$ then $\mathcal{P}=\mathcal{P}_1\star (w,u)\star \overline{\mathcal{P}_2}$ is an $\mathcal{M}$-augmenting path. This follows since $\mathcal{P}_1,\mathcal{P}_2$ are star paths and  $v,y$ are by assumption and definition, respectively, unsaturated vertices.\\
 So, we may assume without loss of generality that $\mathcal{P}_1\cap \mathcal{P}_2 \neq\emptyset$.  Now lemma 7 applies and so we either have a star path from $v$ to $u$ or an $\mathcal{M}$-augmenting path from $v$ to $y$. But the maximality of $\mathcal{M}$ forces the existence of a star path from $v$ to $u$ and that is a contradiction to our assumption that $u\notin T_v$.
Hence $\mathcal {C}$ cannot have more than one unsaturated vertex.\hspace{\stretch{1}} $\square$
 
 We next show an application of theorem 8 to prove Gallai's lemma for Factor Critical graphs. 
\begin{thm} (Gallai's Lemma):  If $G$ is a simple connected graph such that for all $v \in V(G)$, $\nu (G) = \nu (G\setminus {v})$ ($G\setminus {v}$ is the induced subgraph of $G$ on $V(G)\setminus \{v\}$) then $|V(G)|= 2\nu +1$. \end{thm}
 \textbf{Proof}: 
 We first make the following claim: For a connected graph $G$ such that for all $v \in V(G)$, $\nu (G) = \nu (G\setminus {v})$ then $V(G)\subseteq Star(G,\mathcal{M})$ for any maximum matching $\mathcal M$ of $G$.\\
  Suppose not. Then there exists a vertex $v \in V(G)$ such that $v \notin Star(G,\mathcal{M})$. By definition, $\mathcal{M}$-unsaturated vertices are star vertices, hence $v$ is covered by $\mathcal{M}$.  This implies the existence of an edge $e_{vu} \in E(G) \cap \mathcal {M}$. \\
Now consider $G\setminus {v}$. $\mathcal {M} \setminus e_{vu}$ is not a maximum matching of $G\setminus {v}$  as $\nu (G) = \nu (G\setminus {v}) > |\mathcal {M} \setminus e_{vu}|$. So there is an augmenting path in $G\setminus {v}$ for $\mathcal {M} \setminus e_{vu}$. Any such path would be an augmenting path of $G$ for $\mathcal {M}$ as well, unless the path terminates at $u$. Now $u$ is not covered by the matching of $G\setminus {v}$, so, if there is an augmenting path $\mathcal {P}$ in $G\setminus {v}$ for $\mathcal {M} \setminus e_{vu}$ terminating at $u$, there would be a star path $\mathcal {P} \star (u,v)$ in $G$ terminating at $v$. This contradicts $v \notin Star(G,\mathcal{M})$.\\
Now theorem 8 gives us $|V(G)|= 2\nu +1$.\hspace{\stretch{1}} $\square$ 

\begin{prop}: Let $G$ be a simple graph. $\nu (G) = \nu (G\setminus {v})$ for $\forall $ $v \in V(G)$ if and only if $V(G)=Star(G,\mathcal{M})$ for any maximum matching $\mathcal {M}$ of $G$.\end{prop}
\textbf{Proof}: We first show the ``only if" part. Let $\mathcal {M}$ be a maximal matching of $G$ satisfying $V(G)=Star(G,\mathcal{M})$. If a vertex $v \in V(G)$ is not covered by $\mathcal{M}$ then obviously $\nu (G) = \nu (G\setminus {v})$. Now let $u \in V(G) \cap V(\mathcal {M})$. So there is a star path $\mathcal {P} = (v, x_1, x_2, \ldots , w, u)$ in $G$ for some $w \in V(G) \cap V(\mathcal {P})$ with $e_{uw} \in \mathcal {M}$ and for some $v$ not covered by $\mathcal{M}$. Then $(v,x_1,x_2, \ldots , w) $ is an augmenting path in $G\setminus u$ for the matching $\mathcal {M} \setminus  e_{uw}$. That implies $\nu (G \setminus u) > \nu(G) -1$. Hence $\nu(G \setminus u) = \nu(G)$. \\
Now we prove the ``if" part. Recall that a vertex not covered by $\mathcal {M}$ is a star vertex relative to $\mathcal {M}$, so we only need to consider those vertices covered by $\mathcal{M}$. Let $v\in V(G)$ be a vertex covered by $\mathcal{M}$. Then there exists a $w \in V(G)$ such that $e_{wv} \in \mathcal{M}$. Hence there exists an augmenting path $\mathcal {P}$ in $G \setminus v$ for $\mathcal{M} \setminus e_{wv}$ as $\nu(G\setminus v)=\nu(G)$ by the assumption. The path $\mathcal{P}$ contains $w$ since otherwise the path $\mathcal{P}$ would be an augmenting path for $\mathcal{M}$ in $G$, which is impossible. Also note that $w$ is not covered by $\mathcal{M} \setminus e_{wv}$ in $G \setminus v$ so that an augmenting path for $\mathcal{M} \setminus e_{wv}$ can terminate at $w$. So we have a star path $\mathcal{P} \star (w, v)$ in $G$. Hence $v$ is a star vertex of $G$ relative to $\mathcal{M}$.   \hspace{\stretch{1}} $\square$ 
\section{The Transformed graph}
Let $G$ be a graph in $\mathcal{F}(d,m)$. We transform the graph $G$ into a `better-structured' graph, i.e., we obtain another graph whose structure makes the estimation of the number of edges an easier task, which has at least as many edges as $G$ and  which is again a member of $\mathcal{F}(d,m)$. In more precise terms, we seek a graph $G^{\textrm{final}} \in \mathcal{F}(d,m)$ satisfying:
\begin{enumerate}
\item  $\nu(G^{\textrm{final}})=\nu(G)$,
\item  $|E(G^{\textrm{final}})|\geq|E(G)|$,
\item  $ \Delta(G^{\textrm{final}})<d$.  \end{enumerate}
  We go about this task in algorithmic fashion by transforming $G$ in several stages where each intermediate graph is again a member of $\mathcal{F}(d,m)$. But before we describe the process, we set up some more notation.

Let $\mathcal{S}_k$ denote the `claw' $K_{1,k}$ where $K_{1,k}$ denotes the complete bipartite graph with vertex classes of size $1$ and $k$ respectively. For any simple graph $G$, we define \begin{displaymath}\mathcal{C}_{\mathcal{S}_k}(G) :=\{\mathcal{C} \leq G: \mathcal{C}\textrm{\ is a connected component of\ } G, \mathcal{C}\cong K_{1,k}\}.\end{displaymath} For any simple graph $G$ and vertices $u\in V_0(G), v\in V(G)$, (recall that $V_0(G)$ is the set of isolated vertices of $G$) we define another graph $G'$ with $V(G'):=V(G),E(G'):=E(G)\cup e_{vu}$. We denote this graph by $G'=G\oplus e_{vu}$. Note that this is an associative operation (if done in succession for different vertices of $V_0(G)$). We denote the neighborhood of $v\in V(G)$ by $N_G(v)$.\\
In a similar vein, for a simple graph $G$, a subgraph $G'$ of $G$ and $v\in V(G)\cap V(G')$, we denote by $G\ominus E(v,G')$ the graph with vertex set $V(G)$ and edge set $E(G)\setminus \{ e_{vw}| w \in N_{G'}(v)\}$. \\
\textbf{Remark}: Suppose $v$ is a vertex in $G\setminus Star(G, \mathcal{M})$ and $u\in V_0(G)$, then $\nu(G)=\nu(G')$ where $G'=(G\oplus e_{vu})\ominus E(v,G)$.
 Clearly by  Proposition 6, $\nu(G') \leq \nu(G\oplus {e_{vu}})=\nu(G)$. Since $v\notin Star(G,\mathcal{M})$, $v$ is necessarily a saturated vertex relative to $ \mathcal{M}$. Therefore let $e_{vw}$ denote the matching edge incident at $v$. Since we can always define a matching $\mathcal{M}':=\mathcal{M}\setminus\{e_{vw}\}\cup\{e_{vu}\}$ in $G'$, it follows that $\nu(G')\geq \nu(G)$.\\ \\
  Suppose $G\in\mathcal{F}(d,m)$ and let $\mathcal{M}$ be a maximum matching of $G$ then $|\mathcal{M}|=m-1$ by proposition 2. We assume without loss of generality that $|V_0(G)|=2(m-1)(d-1)$. We write $V_0(G)=\displaystyle\bigcup_{v\in V(\mathcal{M})} \mathcal{T}_v$ such that $|\mathcal {T} _v| = d-1$ for all $v\in V(\mathcal{M})$. Note that we necessarily have $\mathcal{T}_u\cap\mathcal{T}_v=\emptyset$.\\ \\
  Let $G_0:=G$ and $\mathcal{M}_0:=\mathcal{M}$. Set 
  \begin{eqnarray}
  G_0'&:=&G_0\setminus\mathcal{C}_{\mathcal{S}_{d-1}}(G_0)\\
  \mathcal{M}_0'&:=& E(G_0')\cap \mathcal{M}_0\\
  t_0&:=&|\mathcal{C}_{\mathcal{S}_{d-1}}(G_0)|.\end{eqnarray}
  
  Note that $|\mathcal{M}_0'|=m-1-t_0$ and $\mathcal{M}_0'$ is a maximal matching of $G_0'$.\\ \\
  If $V(G_0')\setminus Star(G_0',\mathcal{M}_0')=\emptyset$ , then set $G^{\textrm{final}}= G_0$ and we are through.\\ \\
   If not, let $v\in V(G_0')\setminus Star(G_0',\mathcal{M}_0')$ and let $\mathcal{T}_v:=\{w_1,w_2,\ldots,w_{d-1}\}\subseteq V_0(G_0')$.\\
   Now let, $G_1':= G_0'\oplus e_{vw_1}\oplus e_{vw_2}\cdots\oplus e_{vw_{d-1}}$. By Proposition 6, it follows that $\nu(G_1')=\nu (G_0')$. Now let $G_{1*} := G_1'\ominus E(v,G_0')$ . Note that $\nu (G_{1*} )= \nu(G_1')$ by the remark made above. Now set $G_1:=G_{1*}\cup \mathcal{C}_{\mathcal{S}_{d-1}}(G_0)$.\\ \\
   We now describe a maximal matching $\mathcal{M}_1$ of $G_1$:\\
     $\mathcal{M}_1:=\{E(G_1)\cap\mathcal{M}_0'\}\cup \{e_{vw_1}\}\cup\{E(\mathcal{C}_{\mathcal{S}_{d-1}}(G_0))\cap\mathcal{M}_0\}$.\\    
     It is clear that this is a maximum matching of $G_1$ as it gives a maximum matching in every connected component of $G_1$. \\ \\
   Now we iteratively define $G_k$, and a maximal matching $\mathcal{M}_k$ of $G_k$ in a similar manner. Note that $\Delta(G_k)<d$ and that $\nu(G_k)=\nu(G_{k-1})=\cdots=\nu(G_0)$ and hence $G_k\in\mathcal{F}(d,m)$.\\
   This procedure terminates after a finite number of iterations (in fact in at most $m-1$ iterations) since after each iteration the matching size in a component having a vertex which is not a star vertex decreases by one, namely if $G_{k}\setminus\mathcal{C}_{\mathcal{S}_{d-1}}(G_k)$ has a vertex which is not a star vertex, then after an iteration, $\nu((G_{k+1}\setminus\mathcal{C}_{\mathcal{S}_{d-1}}(G_{k+1}))$ decreases by one. The final graph shall be $G^{\textrm {final}}$ with a maximal matching $\mathcal {M^{\textrm {final}}}$. It is clear that $G^{\textrm{final}}$ satisfies conditions $(1),(2)\ \textrm{and}\ (3)$ at the beginning of the section. 
\section{A bound on $|E(G^{\textrm{final}})|$}
By the procedure above, we arrive at a graph $G^{\textrm {final}}$ which has two kinds of connected components:
\begin{enumerate}
\item $K_{1,d-1}$,
\item Components whose vertices are all star vertices (i.e. factor-critical components). Let us denote by $\mathcal {L}$, the set of factor-critical components in $G^{\textrm {final}}$.\end{enumerate}
By theorem 8, we know that for a component $\mathcal {C} \in \mathcal {L},\ \mathcal {C}$ has exactly one unsaturated vertex. Hence if $\mathcal {C}\in\mathcal{L}$ has $r$ matching edges of a maximum matching $\mathcal {M^{\textrm {final}}}$ then,\\
 \begin{displaymath} |V(\mathcal {C})|= 2r+1\textrm{\ and\ } |E(\mathcal {C})|\leq\min\{(2r+1)r,\lfloor {\frac {(2r+1)(d-1)}{2}}\rfloor\}.\end{displaymath}
Thus, if we have $t$ components isomorphic to $K_{1,d-1}$ and $k$ components $\mathcal{C}_i\in\mathcal{L}$, with $|\nu(\mathcal{C}_i)|=r_i$, then we have $|E(G)|\leq e(G):=(d-1)t+\sum_{i=1}^k \min\{(2r_i+1)r_i,\lfloor {\frac {(2r_i+1)(d-1)}{2}}\rfloor\}$ subject to the constraint, $t+\sum_{i=1}^k r_i =m-1$.\\ \\
 Note that $e(G)$ can simply be regarded as a function in the parameters $(t,k,\{r_i\}_{i=1,2\ldots,k})$.\\
 We look to calculate $e_0$, the maximum value of $e(G)$ subject to the linear constraint, $t+\sum_{i=1}^k r_i =m-1$.\\ 
At this juncture, it is worth emphasizing that the rest of this section is an elaboration of a simple idea and can be written more concisely. However, we have included all the minute details for the sake of completeness.\\ 
We start with a few observations in order to rewrite the linear constraint and the expression for the number of edges.
\begin{itemize}
\item  Note that for $r_i\geq\lceil\frac{d-1}{2}\rceil$, we have \begin{displaymath}(2r_i+1)r_i\geq (2r_i+1)\lceil\frac{d-1}{2}\rceil\geq\frac{(2r_i+1)(d-1)}{2}\geq\lfloor\frac{(2r_i+1)(d-1)}{2}\rfloor.\end{displaymath}
 \item If for any $i, r_i> \lceil\frac{d-1}{2}\rceil$, then we can redefine parameters $t$ and $r_i$ in the expression of $e(G)$ and keep the value of $e(G)$ intact. More precisely, let $t'=t+1,r_i=r_i-1$. Since the linear  
constraint is satisfied, the change in $e_0$ is 

\begin{displaymath}
(d-1)(t+1)+\lfloor\frac{(2r_i-1)(d-1)}{2}\rfloor -(d-1)t-\lfloor\frac{(2r_i+1)d-1}{2}\rfloor.
\end{displaymath}\end{itemize}

This simplifies to $ d-1+\lfloor\frac{(2r_i-1)(d-1)}{2}\rfloor-\lfloor\frac{(2r_i+1)d-1}{2}\rfloor=d-1+\lfloor(r_i-1)(d-1)+\frac{d-1}{2}\rfloor-\lfloor r_i(d-1)+\frac{d-1}{2}\rfloor=d-1-(d-1)=0$.\\ \\ Hence we can assume without loss of generality that $r_i\leq \lceil\frac{d-1}{2}\rceil$ for all $1\leq i\leq k$.\\ 
Let $J:=\{1\leq i\leq k\ |r_i=\lceil\frac{d-1}{2}\rceil\}$ and let $I:= \{1,2\ldots,k\}\setminus J$.\\ \\
We consider the cases of $d$ odd and $d$ even.\\ \\
 CASE I: $d$ is odd: Let $d=2j+1$ for some non-negative integer $j$. Then \begin{displaymath} e_0=t(2j)+\sum_{i\in I} (2r_i+1)r_i +(2j+1)j|J|\end{displaymath} with $I,J$ as described above, subject to the constraint,
 \begin{displaymath} t+\sum_{i\in I} r_i+ j|J|=m-1\end{displaymath}
 with $r_i\leq j-1$ for all $i$.\\
  Suppose $r_1>0$. Let $t':=t+r_1, r_1'=0$ and $r_i'=r_i$ for $i>1$. Consider the parameters $(t',\{r_i'\}, J)$. If $e'$ denotes the corresponding value of $e(G)$(as a function of these values), then $e'-e_0=(t+r_1)(2j)-2jt-(2r_1+1)r_1=r_1(2j-2r_1-1)\geq r_1(2j-2j+2-1)>0$ and that contradicts the maximality of $e_0$. \\ 
  Hence $r_i=0\ \forall\  i\in I$.\\ \\
  Thus we have the linear constraint $t+j|J|=m-1$ and we wish to maximize $(2j)t+(2j^2+j)|J|$.\\ It follows from elementary calculus that the maximum occurs at one of the extreme points, i.e., when $|J|=0,t=m-1$ or when $|J|=\lfloor\frac{m-1}{j}\rfloor, t=m-1-j\lfloor\frac{m-1}{j}\rfloor$. It follows that the maximum occurs precisely when $|J|=\lfloor\frac{m-1}{j}\rfloor, t=m-1-j\lfloor\frac{m-1}{j}\rfloor$ and maximum $e_0=2j(m-1)+j\lfloor\frac{m-1}{j}\rfloor$.\\ \\
  CASE II: $d$ is even. Suppose $d=2j$. In this case  \begin{displaymath} e_0=t(2j-1)+\sum_{i\in I} (2r_i+1)r_i +(2j^2-1)|J|\end{displaymath} subject to the linear constraint \begin{displaymath} t+\sum_{i\in I} r_i+ j|J|=m-1,\end{displaymath} again with $r_i\leq j-1$ for all $i$.\\
  As before, suppose $r_1>0$. Then, defining $t',r_i'$ and $e'$ exactly as before, we note that $e'-e_0=(t+r_1)(2j-1)-(2j-1)t-(2r_1+1)r_1=r_1(2j-2r_1-2)\geq 0$. So here we can assume without loss of generality that $r_i=0$ for all $i\in I$(note that equality can occur in the above chain of inequalities).\\ \\
  Once again we are reduced to the case of maximizing $t(2j-1)+(2j^2-1)|J|$ subject to the constraint  $t+ j|J|=m-1$. Exactly as before we see that the maximum occurs at one of the extremities, i.e., at $t=m-1, |J|=0$ or when $t=m-1-j\lfloor\frac{m-1}{j}\rfloor,|J|=\lfloor\frac{m-1}{j}\rfloor$. It is again trivial to see that the maximum occurs when $t=m-1-j\lfloor\frac{m-1}{j}\rfloor,|J|=\lfloor\frac{m-1}{j}\rfloor$ and in this case $e_0=(2j-1)(m-1)+(j-1)\lfloor\frac{m-1}{j}\rfloor$.\\ \\
  So we finally have in either case, $e_0=(d-1)(m-1)+\lfloor\frac{m-1}{\lceil\frac{d-1}{2}\rceil}\rfloor \lfloor\frac{d-1}{2}\rfloor$.  
 \section{$\mathcal{F}(d,m)$-Graphs with $e(d,m)$ edges}
 Having maximized the quantity $e(d,m)$, we now turn our attention back to graphs in the family $\mathcal{F}(d,m)$ that attain this bound. To that end, we first
 discuss construction of a factor critical component $\mathcal {C}$ with $\nu (\mathcal {C})=\lceil\frac{d-1}{2}\rceil $, $\Delta(\mathcal {C})< d$ and $|E(\mathcal{C})|$ as large as possible. For $d=2j+1$ take $\mathcal{C}$ to be a complete graph on $2j+1$ vertices. If $d=2j$ then consider a complete graph $\mathcal {D}$ on $2j$ vertices and then remove $j$ alternate edges of a cycle of $2j$ edges from $\mathcal{D}$. We call this new graph $\mathcal {E}$. Now to obtain $\mathcal{C}$, we connect any of the $2j-1$ vertices of $\mathcal{E}$ to a new, hitherto isolated vertex,  $v \notin V(\mathcal{E})$. We also claim that a factor critical graph $\mathcal {C}$ with $\nu(\mathcal{C})=\lceil\frac{d-1}{2}\rceil $, $\Delta(\mathcal{C})<d$ and $|E(\mathcal{C})|=\lfloor{\frac{(2(\lceil\frac{d-1}{2}\rceil )+1)(d-1)}{2}}\rfloor$ is unique up to isomorphism. If $d$ is odd it is a complete graph on $2(\lceil\frac{d-1}{2}\rceil )+1$ vertices and hence unique. If $d=2j$ for $j>1$ then take $\mathcal{C}$ among the $2j+1$ vertices of $\mathcal{C}$ there is a unique vertex $v$ of degree $2j-2$. Hence there is a vertex $u$ in $V(\mathcal{C})$ which is not a neighbor of $v$. $\mathcal{C} \setminus u$ would be a regular graph of degree $2j-2$ on $2j$ vertices and hence its complement is simply a regular graph of degree one, namely, a matching of a complete graph on $2j$ vertices. This demonstrates the uniqueness of $\mathcal{C}$ and the fact that $\mathcal{C}$ can be obtained from $\mathcal{E}$ as described above.\\
 Now we mention the cases where the extremal graphs (with no isolated vertices) in $\mathcal{F}(d,m)$ are unique. For $d=2$, $G$ is $m-1$ copies of $K_2$.  If $m=2$ and $d\neq 4$ then it is easy to see that a graph with $e(d,2)= d-1$ edges is the complete bipartite graph $K_{1,d-1}$. Suppose $\lceil\frac{d-1}{2}\rceil$ divides $m-1$; in this case, $t=0$. Hence all the vertices of any maximal graph would be star vertices relative to any of its maximum matchings. That implies that $G=G^{\textrm{final}}$ and hence the graph $G$ with maximum number of edges satisfying the conditions $\nu<m$ , $\Delta< d$ would be the unique graph (with no isolated vertices) with $|J|=\lfloor\frac{m-1}{\lceil\frac{d-1}{2}\rceil}\rfloor$ components isomorphic to $\mathcal{C}$, where $\mathcal{C}$ is the graph described in the previous paragraph. However if $\lceil\frac{d-1}{2}\rceil$ does not divide $m-1$, $m>2$ and $d>2$ then many graphs achieve the bound $e(d,m)$. In this case $t \neq 0$ and $G^{\textrm{final}}$ may include at least two components:
 \begin{itemize}
 \item two components isomorphic to $K_{1,d-1}$, 
 \item one component isomorphic to $K_{1,d-1}$ and one component isomorphic to $\mathcal {C}$ as described above. 
 \end{itemize}
 
 We can always coalesce the two components mentioned above into a single component without loss of number of edges. If there are two components $H_1,H_2$ isomorphic to $K_{1,d-1}$, one could remove an edge of $H_1$ and then connect the vertex of degree $d-2$ to any of the degree one vertices of $H_2$. Note that the new graph thus formed has exactly the same number of edges and so is again a member $\mathcal{F}(d,m)$. If in $G^{\textrm {final}}$, there is a component isomorphic to $K_{1,d-1}$ and a component isomorphic to $\mathcal{C}$ then we can form a new component with $2(\lceil\frac{d-1}{2}\rceil +1)+1$ vertices and maximum degree less than $d$, having number of edges equal to the total number of edges in the original two components.\\
 In the special case $d=m=s$, we see that the maximum is $e(s,s)=s(s-1)$ when $s$ is odd and $\lfloor\frac{(2s-1)(s-1)}{2}\rfloor$ when $s$ is even. Furthermore, when $s$ is odd and $s=d=m$ then clearly, $\lceil\frac{d-1}{2}\rceil=\frac{s-1}{2}$ divides $m-1=s-1$ and so in this case the maximum is attained when the graph is isomorphic to two disjoint copies of $K_s$. In the case of $s$ even, there are several graphs that attain the bound.
\section{Acknowledgements}
 We would like to thank Dr. \'Akos Seress for his valuable comments, help with direction of research and guidance in re-writing some parts of the paper.

% \appendix

% \section{}
% \label{}

% Bibliographic references with the natbib package:
% Parenthetical: \citep{Bai92} produces (Bailyn 1992).
% Textual: \citet{Bai95} produces Bailyn et al. (1995).
% An affix and part of a reference:
%   \citep[e.g.][Ch. 2]{Bar76}
%   produces (e.g. Barnes et al. 1976, Ch. 2).

\end{document}